\documentclass[preprint,11pt]{elsart}
\usepackage{amssymb}
\usepackage[T1]{fontenc}
\usepackage[ansinew]{inputenc}

\newtheorem{lemma}[thm]{Lemma}

\newtheorem{definition}[thm]{Definition}
\newcommand{\chro}{\stackrel{\longrightarrow}{{\rm exp}}\int}

\newcommand{\Ad}{\rm Ad}
\newcommand{\ad}{\rm ad \hspace{1pt}}

\newcommand{\spn}{\mbox{Span }}

\begin{document}

\begin{frontmatter}



\title{Controlling Multiparticle System on the Line, II - Periodic case\thanksref{label1}}
\thanks[label1]{This work has been partially supported by MIUR,
Italy, via PRIN 2006019927}

\author{Andrey Sarychev}
\ead{asarychev@unifi.it}
\address{DiMaD, Universit\`a di Firenze, v.C.Lombroso 6/17, Firenze, 50134, Italia}

\begin{abstract}
As in \cite{Sar}  we consider classical system of interacting
particles $\mathcal{P}_1, \ldots , \mathcal{P}_n$ on the line with
only neighboring particles involved in interaction. On the contrast
to [7] now {\it periodic boundary conditions} are imposed onto the
system, i.e. $\mathcal{P}_1$ and $\mathcal{P}_n$ are considered
neighboring. Periodic Toda lattice would be a typical example. We
study  possibility to control periodic multiparticle systems by
means of forces applied to just few of its particles; mainly we
study system controlled by  single force. The free dynamics of
multiparticle systems in periodic and nonperiodic case differ
substantially. We see that also the controlled periodic
multiparticle system does not mimic its non-periodic counterpart.

Main result established  is  global controllability by means of
single controlling force of the multiparticle system with ageneric
potential of interaction. We study the nongeneric potentials for
which controllability and accessibility properties may lack. Results
are formulated and proven in Sections~\ref{mporb},\ref{mpcon}.
\end{abstract}

\end{frontmatter}

\section{Introduction}
\label{intro} \markboth{A.V.Sarychev}{Controllability of
Multiparticle Systems}

Consider classical system of $n$ interacting particles
$\mathcal{P}_1, \ldots , \mathcal{P}_n$ moving on the  line with
only neighboring particles being  involved in the interaction. Let
$q_k$ be the coordinate of the $k$-th particle and $p_k$ - its
momentum. We assume the potential of this interaction to be
\begin{equation}\label{pot}
\Phi(q_1-q_{2})+\Phi(q_2-q_{3})+\cdots
+\Phi(q_{n-1}-q_n)+\Phi(q_{n}-q_1),
\end{equation}
where $\Phi: \mathbb{R} \rightarrow \mathbb{R}$ is real analytic,
bounded below function
\begin{equation}\label{grof}
\lim_{y \rightarrow +\infty}\Phi(y)=+\infty.
\end{equation}

The difference with  nonperiodic case, studied in \cite{Sar}, is due
to the presence of the last addend in (\ref{pot}) which accounts for
neighboring of $\mathcal{P}_1$ and $\mathcal{P}_n$. This extra
addend leads to  a substantial change of dynamics. For example in
celebre and extensively studied case of Toda lattice, in which
interaction potential equal to $\Phi(x)=e^{2x}$, the distances
between particles are known (\cite{Mos}) to  tend to infinity in
nonperiodic case, while in periodic case the particles are involved
in quasiperiodic motion on compact isoenergetic surfaces. Below  we
will see that the controlled dynamics in the periodic case also
differs from nonperiodic controlled dynamics studied in Part I.

The dynamics of multiparticle system with the  potential (\ref{pot})
is described by   Hamiltonian system of equations with the
Hamiltonian
\begin{equation}\label{ham-fi}
H(q,p)=\frac{1}{2}\sum_{k=1}^np_k^2+\sum_{j=1}^{n}\Phi(q_j-q_{j+1}).
\end{equation}
Corresponding equations are
\begin{eqnarray}
\dot{q}_k=\frac{\partial H}{\partial p_k}=p_k, \ k=1, \ldots ,n, \label{dqq}\\
\dot{p}_k=-\frac{\partial H}{\partial
q_k}=\phi(q_{k-1}-q_k)-\phi(q_{k}-q_{k+1}), \ k=2, \ldots , n, \label{dpq}\\
\dot{p}_1=-\frac{\partial H}{\partial
q_1}=\phi(q_{n}-q_1)-\phi(q_{1}-q_{2}). \label{dp1}
\end{eqnarray}
In (\ref{dpq}),(\ref{dp1}) and further on  $\phi=\Phi '$ is the
derivative of $\Phi$. Besides for unification of notation  we assume
in (\ref{ham-fi}), (\ref{dpq}) and (\ref{dp1})
$$ q_0=q_n,\  q_{n+1}=q_1.$$

The control will be realized by  a force, which we choose to act on
the particle  $\mathcal{P}_1$.  In the presence of the control the
equation (\ref{dp1})  becomes
\begin{equation}\label{conp1}
\dot{p}_1=\phi(q_{n}-q_1)-\phi(q_{1}-q_{2})+u(t),
\end{equation}
where $u(\cdot)$ stays for the controlling force. The equations
(\ref{dqq}),(\ref{dpq}) remain unchanged. We call the model {\it
single forced periodic multiparticle system}.

We wish to study controllability properties   of the controlled
multiparticle periodic system (\ref{dqq})-(\ref{dpq})-(\ref{conp1}).

\begin{definition}
\label{gcdef} System (\ref{dqq})-(\ref{dpq})-(\ref{conp1}) is
globally controllable if for given pair of points
$\tilde{x}=(\tilde{q},\tilde{p}), \hat{x}=(\hat{q},\hat{p})$ of its
state space there exists an admissible (measurable essentially
bounded) control $u(\cdot)$ which steers the system from $\tilde{x}$
to $\hat{x}$ in time $\theta>0. \ \square$
\end{definition}

The controlled multiparticle system
(\ref{dqq})-(\ref{dpq})-(\ref{conp1}) is a particular case of
control affine system of the form

\begin{equation}\label{pcaf}
    \frac{dx}{dt}=f(x)+g(x)u,
\end{equation}
where  the {\it controlled vector field}
 $g$ and   the uncontrolled vector field
$f$ - the {\it drift} - are defined as
\begin{equation}\label{pervef}
g=\partial / \partial p_1, \
f=\sum_{k=1}^{n}p_k\frac{\partial}{\partial q_k}+\sum_{k=1}^{n}
\left(\phi(q_{k-1}-q_k)-\phi(q_{k}-q_{k+1})\right)\frac{\partial}{\partial
p_k}.
\end{equation}

In Part I we observed  that in non periodic case global
controllability is in general non-achievable by means of one
controlling force and is achievable by means of two controlling
forces applied to the "extreme" particles $\mathcal{P}_1$ and
$\mathcal{P}_n$. It is immediate to conclude (see
Subsection~\ref{lrper}) that the periodic multiparticle system is
also globally controllable by means of two forces.

We are going to prove instead that {\small\it in the periodic case
for a generic potential $\Phi$  global controllability is achievable
by means of single controlling force} (Theorem~\ref{pecon} in
Section~\ref{mpcon}). The proof is split into two parts. First
 we establish full dimensionality of the {\it orbit} of a single forced
 multiparticle system.

 An orbit $\mathcal{O}_{\tilde{x}}$
 of control system is the minimal invariant manifold for the
 control system, whenever one starts from the initial point
 $\tilde{x}$ and proceeds with controlled motion in direct (positive) and reverse (negative) time.
 We will prove  that for a
  {\it generic potential} $\Phi$ the orbits of the control system (\ref{pcaf})-(\ref{pervef})
  coincide with the state space $\mathbb{R}^{2n}$. This is
   done  (Subsection~\ref{lrper})  by verification of
     {\it bracket generating property} of the couple of vector fields
   $\{f,g\}$.
 This property may fail for some potentials;
 in  Subsection~\ref{ldpr} we provide an example of
 {\it low-dimensional orbits for a specific potential} $\Phi$. In
 Subsection~\ref{ldnp} we return for a moment  to
 nonperiodic case and provide an example of {\it low-dimensional orbit
 for
 nonperiodic} system whenever controlling force is applied to a particle $\mathcal{P}_j$ with $j
 \neq 1,j \neq n$.

Once  full dimensionality of an orbit is established, one has to
deal with another difficulty. Positive invariant set of a control
system ({\it attainable set}) is often a proper subset of the
respective orbit.
  The reason for this is actuation of the drift vector field $f$, which may drive
  the system in certain direction without a possibility to compensate this drift
  by action of any control.

  In some exceptional
 cases such compensation is possible.
One of these cases is represented by Bonnard-Lobry theorem
 (\cite{Bo}), whose main assumption is
{\it recurrence property of dynamics of the noncontrolled motion}.

 In the nonperiodic case, treated in \cite{Sar}, we arranged a simple design of  feedback controls which
 modified
 the noncontrolled dynamics  in such a way that all  its trajectories  became recurrent.
    Such design was only possible with  two controls available.

  In the periodic case we get instead a  property of {\it constrained recurrence} for
   the dynamics of non-controlled motion: the  dynamics is recurrent on a
   hyperplane
 of zero momentum $\Pi: \ p_1+ \cdots +p_n=0. $
The hyperplane is invariant with respect to free dynamics, but is
not invariant with respect to  controlled dynamics. Therefore one
can not remain in $\Pi$, whenever nonzero control is employed, and
we can not use the recurrence property when one is  outside $\Pi$.
We will adapt the technique of {\it Lie extensions} for overcoming
this difficulty  and establishing global controllability.

\section{Orbits and accessibility property for single forced multiparticle system}
\label{mporb}

We study single-forced periodic multiparticle system, or, the same,
control-affine system (\ref{pcaf})-(\ref{pervef})in  the {\it state
space} $\mathbb{R}^{2n}$.

We start with computation in the next subsection  of the orbits of
this control-affine system. Recall that one obtains orbit
$\mathcal{O}_{\tilde{x}}$  by taking vector fields $f^{u^j}=f+u^j g$
with $u^j \in \mathbb{R}$ constant, and acting on $\tilde{x} \in
\mathbb{R}^{2n}$ by the compositions
\begin{equation}\label{calP}
P=  e^{t_1f^{u^{j_1}}} \circ \cdots \circ e^{t_Nf^{u^{j_N}}}, \ t_1,
\ldots , t_N \in \mathbb{R},
\end{equation}
where $e^{tX}$ stays for the flow of the vector field $X$.

 According to Nagano theorem (\cite{ASkv,Jur} an orbit
$\mathcal{O}_{\tilde{x}}$ is an immersed submanifold of
$\mathbb{R}^{2n}$ and the tangent space to this manifold at a point
$x \in \mathcal{O}_{\tilde{x}}$ is obtained by evaluation at $x$ of
the vector fields from the Lie algebra $\mbox{Lie}\{f,g\}$ generated
by $f$ and $g$.

\begin{definition}
\label{bgp} A family $\mathcal{F}$ of vector fields is called
bracket generating at point $x \in \mathbb{R}^{2n}$ if the
evaluation at $x$ of the vector fields from $\mbox{Lie}\{F\}$
coincides with $\mathbb{R}^{2n}. \ \square$
\end{definition}

An {\it attainable set} $\mathcal{A}_{\tilde{x}}$ of the system
(\ref{pcaf}) from $\tilde{x}$  is the set of points to which the
system can be steered from $\tilde{x}$ by means of an admissible
(measurable, bounded) control. If we require in addition the
transfer time to be equal, or respectively, $\leq$ than  $T$, then
we obtain time-$T$ (respectively time-$\leq T$) attainable set
$\mathcal{A}^T_{\tilde{x}}$ (respectively $\mathcal{A}^{\leq
T}_{\tilde{x}}$). Obviously $\mathcal{A}^T_{\tilde{x}} \subset
\mathcal{A}^{\leq T}_{\tilde{x}} \subset   \mathcal{A}_{\tilde{x}}
$; also $\mathcal{A}_{\tilde{x}}$ is contained in the orbit
$\mathcal{O}_{\tilde{x}}$.

If one employs piecewise-constant controls, i.e. takes only positive
times $t_j >0$ in the compositions (\ref{calP}), then one gets  {\it
positive orbit} $\mathcal{O}^+_{\tilde{x}}$ of the system. In
general it is proper subset of $\mathcal{O}_{\tilde{x}}$ and is far
from being a manifold. Obviously $\mathcal{O}^+_{\tilde{x}} \subset
\mathcal{A}_{\tilde{x}}$.

\begin{rem}
\label{dense}  It is known  from A.J.Krener theorem
(\cite{ASkv,Jur}), that for each point $\tilde{x}$  positive orbit
$\mathcal{O}^+_{\tilde{x}}$ possesses nonvoid relative interior in
the orbit $\mathcal{O}_{\tilde{x}}$, and moreover
$\mathcal{O}^+_{\tilde{x}}$ is contained in the closure of its
relative interior. A consequence of this theorem is the useful fact
(see \cite{ASkv}) that {\it density} of $\mathcal{O}^+_{\tilde{x}}$
in the orbit $\mathcal{O}_{\tilde{x}}$ implies the coincidence of
$\mathcal{O}^+_{\tilde{x}}$ and $\mathcal{A}_{\tilde{x}}$ with
$\mathcal{O}_{\tilde{x}}. \ \square$
\end{rem}

\subsection{Orbits of single forced periodic multiparticle system}
\label{lrper}

In \cite{Sar} we proved for single forced non periodic multiparticle
system that all the orbits coincide with $\mathbb{R}^{2n}$. In the
periodic case this holds for generic potentials. We prove this fact
in the present Subsection and provide counterexamples in
Subsections~\ref{ldpr},~\ref{ldnp}.

\begin{thm}
\label{clieper} For a generic potential $\Phi$ the system of vector
fields $\{f,g\}$ is bracket generating at each point of the state
space $\mathbb{R}^{2n}$; therefore $\forall \tilde{x} \in
\mathbb{R}^{2n}$ the orbit $\mathcal{O}_{\tilde{x}}$ of single
forced multiparticle periodic system (\ref{pcaf}),(\ref{pervef})
through $\tilde{x}$ coincides with $\mathbb{R}^{2n}. \ \Box$
\end{thm}

The proof is structured in two Lemmas, first of which  mimics
similar result for double forced {\it non periodic} multiparticle
system.

Assume for the moment that  periodic multiparticle system is
controlled by {\em two} forces applied to the particles
$\mathcal{P}_1$ and $\mathcal{P}_n$, i.e. we gain an additional
controlled vector field $g^n=\frac{\partial}{\partial p_n}$. The
additional controlling force appears in  the equation (\ref{dpq})
indexed by $k=n$, which now will take form
\begin{equation}\label{conpn}
\dot{p}_n=\phi(q_{n-1}-q_n)-\phi(q_{n}-q_{1})+v(t).
\end{equation}

\begin{lemma}\label{per2in}
The family of vector fields $\{f,g,g^n\}$ is bracket generating at
each point of $\mathbb{R}^{2n}. \ \square$
\end{lemma}

{\it Proof.} The conclusion of the Lemma follows from the fact that
feedback transformation
$$u \mapsto - \phi(q_{n}-q_{1})+\tilde{u}, \ v \mapsto
\phi(q_{n}-q_{1})+\tilde{v},$$  transforms the equations
(\ref{conp1}),(\ref{conpn}) into respective equations of double
forced {\it nonperiodic} multiparticle system,  whose orbits
coincide with $\mathbb{R}^{2n}$ by results of \cite{Sar}.  This
transformation does not affect bracket generating property, hence
the vector fields $\{f,g,g^n\}$ form a  bracket generating system.
$\square$

Recall the notation: for a vector field $X$ operator $\ad X$ acts on
another vector field $Y$ as $\ad X Y=[X,Y]$. The proof of
Theorem~\ref{clieper} would be  accomplished by  the following
lemma.

\begin{lemma}\label{getpn}
For a generic potential $\Phi$, for each point $\tilde{x} \in
\mathbb{R}^{2n}$:
\begin{equation}\label{incl}
\!\!\!\!\! \mbox{Span}\left\{f(x),g(x),\ad^2f g(x),[\ad f g,\ad^2f
g](x)\right\} \supset \mbox{Span}\left\{f(x),
\frac{\partial}{\partial p_1},\frac{\partial}{\partial p_n}\right\},
\end{equation}
 for all $x$ of an open dense subset of the orbit
$\mathcal{O}_{\tilde{x}}$ of the system (\ref{pcaf}),(\ref{pervef}).
$\square$
\end{lemma}

For those potentials $\Phi$, for which the conclusion of
Lemma~\ref{getpn} is valid, one easily gets the statement of
Theorem~\ref{clieper} proven. Indeed since the system  $\{f,
\frac{\partial}{\partial p_1},\frac{\partial}{\partial p_n}\}$ is
bracket generating at each point, then by Lemma~\ref{getpn} the
system of vector fields $\{f,g\}$ is bracket generating at some
point of each orbit $\mathcal{O}_{\tilde{x}}$. The dimension $\dim
\mbox{Lie}_x\{f,g\}$ of the evaluation at  $x$ of the  Lie algebra
 $\mbox{Lie}\{f,g\}$ is known to be constant along $\mathcal{O}_{\tilde{x}}$ (see
\cite{ASkv,Jur}). Hence we conclude that $\{f,g\}$ is bracket
generating at each point of $\mathbb{R}^{2n}$ and all orbits of the
system (\ref{pcaf}),(\ref{pervef}) coincide with the state space
$\mathbb{R}^{2n}$.

{\sf Proof of Lemma~\ref{getpn}.}  By direct computation $[f,g]=\ad
fg=-\frac{\partial}{\partial q_1}$. Computing the iterated Lie
brackets $\ad^2f g,[\ad^2f g, \ad fg]$ we get
\begin{eqnarray}\label{ad2fg}
\!\!\!\!\!\ad^2f g=\phi'(q_1-q_2)\left(\frac{\partial}{\partial
p_2}-\frac{\partial}{\partial p_1}\right)+ \phi'(q_n-q_1)
\left(\frac{\partial}{\partial p_n}-\frac{\partial}{\partial
p_1}\right), \\
\label{ad2adf}
\!\!\!\!\![\ad^2f g, \ad fg]=\phi''
(q_1-q_2)\left(\frac{\partial}{\partial
p_2}-\frac{\partial}{\partial p_1}\right)- \phi''(q_n-q_1)
\left(\frac{\partial}{\partial p_n}-\frac{\partial}{\partial
p_1}\right).
\end{eqnarray}

 We would  arrive to the
needed conclusion  at each point $x \in \mathcal{O}_{\tilde{x}}$
where  the determinant
\begin{eqnarray}
 det\left(
        \begin{array}{cc}
          \phi'(q_1-q_2) &  \phi'(q_n-q_1) \\
          \phi''(q_1-q_2) & -\phi''(q_n-q_1) \\
        \end{array}
      \right)= -\phi'(q_1-q_2)\phi''(q_n-q_1)- \nonumber \\
      -\phi'(q_n-q_1)\phi''(q_1-q_2)=(\phi'(q_1-q_2))^2\frac{\partial}{\partial
q_1}\frac{\phi'(q_n-q_1)}{\phi'(q_1-q_2)} \label{det}
\end{eqnarray}
is nonvanishing.

 As far as the vector field
$[g,f]=\frac{\partial}{\partial q_1}$ is tangent to any orbit
$\mathcal{O}_{\tilde{x}}$ of (\ref{pcaf}) then we get the conclusion
of the Lemma whenever the determinant (\ref{det}) (an analytic
function) {\it does not vanish identically with respect to $q_1$}.
The latter possibility occurs only if the relation
\begin{equation}\label{phic}
\phi'(q-q_2)=c\phi'(q_n-q)
\end{equation} holds identically with
respect to $q$ (by which we substituted $q_1$) with $c$  constant.
Seeing now  $q_2,q_n,c$ as parameters, we treat  (\ref{phic}) as a
functional equation.

Substituting $q=\frac{q_n+q_2}{2}-t$ into (\ref{phic}) we obtain the
relation

$$\forall t: \ \phi'\left(\frac{q_n-q_2}{2}-t\right)=c
\phi'\left(\frac{q_n-q_2}{2}+t\right),$$ wherefrom $c=\pm 1$.  Then
\begin{equation}\label{evod}
\!\!\!c=\pm 1, \phi'(t)=f\left(t-b\right), \ f \ \mbox{is even or
odd, according to the sign of}\ c,
\end{equation}
and  $b=\frac{q_n-q_2}{2}$.

 (See \cite{Pol} for an alternative
description of solution of (\ref{phic}).)

For  {\it generic} potentials $\Phi$, namely for those whose second
derivative $\Phi''=\phi'$  is {\it not} of the form (\ref{evod}),
the inclusion (\ref{incl}) holds on an open dense subset of any
orbit and then these orbits coincide with the state space
$\mathbb{R}^{2n}. \ \square$

It is interesting to know whether there exist potentials $\Phi$, for
which the system (\ref{pcaf}),(\ref{pervef}) possesses
low-dimensional orbits.  In the next two Subsections we provide such
 examples.

\subsection{Low-dimensional orbits of single-forced periodic multiparticle system}
\label{ldpr}
 Consider a trimer - periodic three-particle system
with the Hamiltonian
$$H=\frac{1}{2}\left(p_1^2+p_2^2+p_3^2\right)+\Phi(q_1-q_2)+\Phi(q_2-q_3)+\Phi(q_3-q_1).$$
and the {\it controlled} dynamics
\begin{eqnarray}\label{trimer}
\dot{q}_i=p_i, \ i=1,2,3; \nonumber \\
\dot{p}_1=\phi(q_3-q_1)-\phi(q_1-q_2)+u, \\
\dot{p}_2=\phi(q_1-q_2)-\phi(q_2-q_3), \
\dot{p}_3=\phi(q_2-q_3)-\phi(q_3-q_1), \nonumber
\end{eqnarray}
where $\phi(q)=\Phi'(q)$.  We assume the derivative $\phi'(q)$ to be
of the form (\ref{evod}) with $f$  {\it even} \footnote{One can
prove that whenever $c=-1$ in (\ref{phic}) and respectively $f$ is
odd in (\ref{evod}) the orbits coincide with $\mathbb{R}^{2n}$}.
Then $\phi(t)=F\left(t-b\right)+c$, where $F(t)$ is a primitive of
even function $f(t)$, and therefore can be chosen an {\it odd
function}. In this case $\phi(t)=F\left(t-b\right)+\phi(b)$.

From the differential equations for $q_2,p_2,q_3,p_3$ in
(\ref{trimer}) we derive
\begin{eqnarray*}
  \frac{d}{dt}\left(q_3-q_2\right)=\left(p_3-p_2\right),  \\
\frac{d}{dt}\left(p_3-p_2\right)=2\phi(q_2-q_3)-
\phi(q_1-q_2)-\phi(q_3-q_1)=
\\ =2F\left(q_2-q_3-b\right)-F\left(q_1-q_2-b\right)-
F\left(q_3-q_1-b\right).
\end{eqnarray*}
Assuming in addition $F(-3b)=0$, or equivalently
$\phi(b)=\phi(-2b)$, we check immediately that the $4$-{\it
dimensional plane}
$$\Pi_b: \ p_3-p_2=0, \ q_3-q_2=2b , $$
is an invariant manifold for the control system  (\ref{trimer}).
Indeed, along  $\Pi_b$
\begin{eqnarray*}
  F\left(q_2-q_3-b\right)=F(-3b)=0,   \\
  F\left(q_1-q_2-b\right)+
F\left(q_3-q_1-b\right)=F\left(q_1-q_2-b\right)+
F\left(q_2-q_1+b\right)=0.
\end{eqnarray*}
 Hence $\forall \tilde{x} \in \Pi_b$ the orbit $\mathcal{O}_{\tilde{x}}$ of the control system
(\ref{trimer}) is contained in $\Pi_b$.

\subsection{Low-dimensional orbits of nonperiodic multiparticle
system} \label{ldnp}

In Part I we mentioned that non periodic single forced multiparticle
system may possess low-dimensional orbits for some  choices of
$\phi$, whenever the control force is applied to a particle
different from $\mathcal{P}_1, \mathcal{P}_n$. Here we provide such
example obtained by a variation on the example of the previous
subsection.

For a {\it non periodic trimer} with the controlling force acting on
the particle $\mathcal{P}_2$ the dynamic equations are
\begin{eqnarray}\label{nptrim}
\dot{q}_i=p_i, \ i=1,2,3;  \\
\dot{p}_1=-\phi(q_1-q_2), \
\dot{p}_2=\phi(q_1-q_2)-\phi(q_2-q_3)+u, \ \dot{p}_3=\phi(q_2-q_3).
\nonumber
\end{eqnarray}
Then
\begin{equation}\label{p31}
\frac{d}{d t}\left(q_3-q_1\right)=p_3-p_1, \ \frac{d}{d
t}\left(p_3-p_1\right)=\phi(q_2-q_3)+\phi(q_1-q_2).
\end{equation}

Let us choose the function $\phi(t)=f(t-b), \  f$ - odd function.
Considering the $4$-dimensional plane $\Pi'_b: \ q_3-q_1=-2b, \
p_3-p_1=0. $ we claim that $\Pi'_b$ is invariant  for the control
system (\ref{nptrim}). Indeed restricting the second one of the
equations (\ref{p31}) to $\Pi'_b$ we conclude
$$\frac{d}{d
t}\left(p_3-p_1\right)=\phi(q_2-q_3)+\phi(q_3+2b-q_2)=f(q_2-q_3-b)+f(q_3-q_2+b)=0,$$
 independently of a choice of control $u(\cdot)$. The first one of the equations (\ref{p31}) restricted
 to $\Pi'_b$ implies: $\frac{d}{d
t}\left(q_3-q_1\right)=\left(p_3-p_1\right)=0$. Thus the
$4$-dimensional plane $\Pi'_b$ contains the  orbits $\mathcal{O}_x$
of the control system (\ref{nptrim}) for each $x \in \Pi'_b$.

\section{Controllability of periodic multiparticle
system by means of a single force} \label{mpcon}

In \cite{Sar} we designed special feedback controls which imposed
{\it recurrent behavior} on  dynamics  of {\it nonperiodic double
forced} multiparticle system.  This  allowed us to apply
Bonnard-Lobry theorem (\cite{Bo}) for proving  global
controllability. The same procedure can be repeated for {\it double
forced periodic case}.

\begin{prop}
 Periodic multiparticle  system  is globally
controllable  by means of controlling forces  applied to the
particles $\mathcal{P}_1,\mathcal{P}_n$. $\Box$
\end{prop}

We are aiming though at a stronger result.

\begin{thm}\label{pecon}
Periodic multiparticle system  with generic interaction potential
$\Phi$ is globally controllable by means of a single force. $\Box$
\end{thm}

\begin{rem}
Genericity assumption for the potential $\Phi$ is the same one,
which appeared in Subsection~\ref{lrper} in the course of
computation of orbits. $\square$
\end{rem}

\begin{rem}
There are no a priori constraints imposed on the magnitude of the
controlling force in the formulation of Theorem~\ref{pecon}.
$\square$
\end{rem}

In the rest of this contribution we {\sf prove Theorem~\ref{pecon}}.

\subsection{Lie extensions}

The following definition is  slight modification of the notion of
Lie saturation introduced by V.Jurdjevic (\cite{Jur})
\begin{definition}Let
$\mathcal{F}$ be a family of  analytic vector fields, and
$\mbox{Lie}(\mathcal{F})$ be the Lie algebra generated by
$\mathcal{F}$. Lie extension $\hat{\mathcal{F}}$  of $\mathcal{F}$
is a family $\hat{\mathcal{F}} \subseteq \mbox{Lie}(\mathcal{F})$
such that
\begin{equation}\label{stsat}
\mbox{clos}{\mathcal A}_{\hat{\mathcal F}}(\hat{x}) \subseteq
\mbox{clos}{\mathcal A}_{\mathcal F}(\hat{x}).
\end{equation}
Any vector field from a Lie extension is called {\it compatible}
with $\mathcal{F}. \ \square$
\end{definition}

We specify some types of Lie extensions.

\begin{prop}\label{clos}
A closure  $clos(\mathcal F)$ of $\mathcal F$ in the Whitney
$C^\infty$-topology is a Lie extension. $\Box$
\end{prop}

This assertion follows from  classical result on continuous
dependence of the solutions of ODE on initial data and the r.-h.
side.

An important kind of extension which underlies theory of relaxed or
sliding mode controls is introduced by the following

\begin{prop}\label{conh}
For a control system $\mathcal F$ its  conic hull
$$\mbox{cone}(\mathcal F)=
\left\{\sum_{j=1}^N \alpha_j f^j | \ \alpha_j \in C^\infty (\mathbb{R}^n), f^j \in \mathcal F,
\ N \in \mathbb{N}, \   \alpha_j \geq 0, \ j=1, \ldots ,
N\right\},$$ is a Lie extension.  $\Box$
\end{prop}

To introduce another  type of Lie extension we define {\it
normalizer} of $\mathcal{F}$.

\begin{definition}[see \cite{Jur}]\label{norm}
Diffeomorphism $P$ is a normalizer for the family $\mathcal F$ of
vector fields if $\forall \hat{x}$:
$$P\left({\mathcal A}_{\mathcal F}(P^{-1}(\hat{x}))\right)
\subseteq \mbox{clos}{\mathcal A}_{\mathcal F}(\hat{x}).  \ \Box
$$

\end{definition}

The following sufficient criterion  is useful for finding
normalizers.

\begin{prop}[\cite{Jur}]\label{invnorm}
Diffeomorphism $P$ is a normalizer for the family $\mathcal F$ if
both $P(\hat{x})$ and $P^{-1}(\hat{x})$ belong to
$\mbox{clos}\left({\mathcal A}_{\mathcal F}(\hat{x})\right), \
\forall \hat{x}. \ \Box$
\end{prop}

Now we define an extension. Recall that {\it adjoint action of
diffeomorphism $P$ on a vector field $f$} results in another vector
field  defined as
$$\Ad Pf(x)=\left.P^{-1}_*\right|_{P(x)}f(P(x)). $$

\begin{prop}
The set
$$\tilde{\mathcal{F}}=\{\Ad P f| \ f \in \mathcal{F}, P \ \mbox{- normalizer of} \ \mathcal{F}\}$$
is a Lie extension of $\mathcal{F}. \ \Box$
\end{prop}

For control-affine of the form (\ref{pcaf}) the family of vector
fields, which determines  polidynamics of such system, is
$\mathcal{F}=\{f+gu| \ u \in \mathbb{R}\}$.

According to Propositions~\ref{clos},\ref{conh} the vector fields
\begin{equation}\label{pmg}
\pm g=\lim_{\theta \rightarrow 0}\theta^{-1}(f+g(\pm \theta )) ,
\end{equation}
are contained in  the closure of the conic hull of $\mathcal{F}$ and
therefore are compatible with $\mathcal{F}$.

By Proposition~\ref{invnorm} each diffeomorphism $e^{ug}$ is a
normalizer of $\mathcal{F}$ and hence there holds thew following

\begin{lemma}
\label{lex} The vector fields $\{\Ad e^{\pm u g}f| \ u \in
\mathbb{R} \}$ are compatible with the control system (\ref{pcaf}).
$\square$
\end{lemma}

\subsection{Lie extension for single forced periodic multiparticle system}
\label{unb}

We will employ Lie extensions  for proving Theorem~\ref{pecon}.

Direct computation of  $e^{u\ad g}f$ for the vector fields
(\ref{pervef}) results in
$$b_u=e^{u\ad g}f=f+u[g,f]; $$
it suffices to note that  $\ad^2g f=[g,[g,f]]=0$.

Consider  vector field the $-f$ and join  it to the vector fields $
b_1, b_{-1}$. The three vector fields are contained in
2-distribution $\mathcal{D}$ spanned by $f$ and $[g,f]$.

Above we introduced the plane of zero momentum $\Pi: \ P=p_1 +
\cdots + p_n=0$, which is invariant under the free motion.

\begin{lemma}
\label{liep0} The hyperplane $\Pi$ is invariant for 2-distribution
$\mathcal{D}$, which is bracket generating on $\Pi. \ \Box$
\end{lemma}

{\small\bf Proof of Lemma~\ref{liep0}.} By direct computation (see
formulae (\ref{ad2fg}),(\ref{ad2adf})) one checks that distribution
$\mathcal{D}$ is tangent to $\Pi$.   Also for {\it generic}
$\phi=\Phi'$,
 (see Subsection~\ref{lrper}):
$$\spn\{\ad^2fg, [\ad^2 fg, \ad fg]\}=\spn\left\{Y^2,Y^n\right\},$$
where $Y^2=\frac{\partial}{\partial p_2}-\frac{\partial}{\partial
p_1},\ Y^n=\frac{\partial}{\partial p_n}-\frac{\partial}{\partial
p_1}$. Again by direct computation
\begin{eqnarray*}
[Y^2,f]=Z^2=\frac{\partial}{\partial q_2}-\frac{\partial}{\partial
q_1}, \  [Y^n,f]=Z^n=\frac{\partial}{\partial
q_n}-\frac{\partial}{\partial q_1},
 \\ %
 \ [Z^2,f]=\frac{\partial}{\partial p_3}-\frac{\partial}{\partial
p_1}, \ [Z^n,f]=\frac{\partial}{\partial
p_{n-1}}-\frac{\partial}{\partial p_1} \ (\mbox{mod }
\spn\{Y^2,Y^n\}).
\end{eqnarray*}

We can arrive to the conclusion of the Lemma by induction. $\Box$

The conic hull of the triple of vector fields $\{-f, b_1, b_{-1}\}$
coincides with $\mathcal{D}$. Hence by Rashevsky-Chow theorem
(\cite{ASkv,Jur})
 for each $\tilde{x} \in \Pi$ positive orbit
$\mathcal{O}^+_{\tilde{x}}$ of this triple  is dense in the orbit of
$\mathcal{D}$,
 equal to $\Pi$. By Remark~\ref{dense} it must coincide with
 $\Pi $.

 Note that
 $$e^{t\Ad (e^{\pm g})f}=\Ad (e^{\pm g})e^{tf}=e^{\pm g} \circ e^{tf}\circ e^{\mp g}.$$
  According to the aforesaid each point of
$\Pi$ is attainable from another point of $\Pi$ by means of
composition of diffeomorphisms from the family
\begin{equation}\label{trifam}
\{e^g \circ e^{tf}\circ e^{-g}, \ e^{-g} \circ e^{tf}\circ e^{g}, \
e^{-tf}, \ t \geq 0\};
\end{equation}
$\Pi$ is invariant under the action of diffeomorphisms
(\ref{trifam}).

 We wish to achieve global controllability on $\Pi$ without
having recourse to $e^{-tf}$.

\begin{prop}
\label{bezpm} Each point of $\Pi$ is attainable from another point
of $\Pi$ by means of compositions of diffeomorphisms from the family
\begin{equation}\label{corfam}
\{e^g \circ e^{tf}\circ e^{-g}, \ e^{-g} \circ e^{tf}\circ e^{g}, \
e^{tf}, \ t \geq 0\}. \Box
\end{equation}
\end{prop}

The proof of the Proposition~\ref{bezpm}, postponed to
Subsection~\ref{propost}, follows the line of the proof of
Bonnard-Lobry theorem (see \cite{Bo,ASkv}) and is based on the
recurrence property of the free motion of the multiparticle system
in the plane $\Pi$. Meanwhile taking it conclusion for granted we
accomplish the proof of global controllability.

\subsection{Proof of global controllability}
\label{gcsfms}

By direct computation one checks that for controlled motion the
total momentum $P$ varies according to the equation $\dot{P}=u(t)$.
Taking two points $(\tilde{q},\tilde{p}), (\bar{q},\bar{p})$ in the
state space,  we can steer, say in time $1$, the system (\ref{pcaf})
from $(\tilde{q},\tilde{p})$ to some point
$(\tilde{q}^0,\tilde{p}^0)$ of $\Pi$ by application of a constant
control $\tilde{u}$. Considering the reverse time dynamics
$\dot{P}=-u(t)$  one ensures the possibility to steer the system
(\ref{pcaf}) in time $-1$ from the point $(\bar{q},\bar{p})$ to a
point $(\bar{q}^0,\bar{p}^0)$ of the plane $\Pi$ by means of another
constant control $\bar{u}$. In direct time the system (\ref{pcaf})
would shift in time $1$  from $(\bar{q}^0,\bar{p}^0)$ to
$(\bar{q},\bar{p})$ under the action of $\bar{u}$.

According to Proposition~\ref{bezpm} one can steer the point
$(\tilde{q}^0,\tilde{p}^0)$ to the point $(\bar{q}^0,\bar{p}^0)$ by
a composition of diffeomorphisms of the form $e^{g},e^{-g},e^{t f},
\ t \geq 0$. Then this composition of diffeomorphisms preceded by
time-$1$ action of the control $\tilde{u}$ and succeeded by time-$1$
action of the control $\bar{u}$ steers the system from
$(\tilde{q},\tilde{p})$ to $(\bar{q},\bar{p})$ in the state space.

According to the limit relation (\ref{pmg}) we can approximate
arbitrarily well the diffeomorphisms  $e^{\pm g}$ in the
composition, we have just described,  by diffeomorphisms
$e^{\theta^{-1}(f \pm g\theta)}$ with sufficiently large $\theta
>0$; these latter are elements of admissible flows $e^{t(f \pm
g\theta)}$. Hence one can steer the point $(\tilde{q},\tilde{p})$ by
an admissible control to a point $(\bar{q}',\bar{p}')$ which is
arbitrarily close to $(\bar{q},\bar{p})$. As far as
$(\bar{q},\bar{p})$ is arbitrarily chosen, one concludes that the
attainable set of the system (\ref{pcaf})-(\ref{pervef}) from each
point $(\tilde{q},\tilde{p}) \in \mathbb{R}^{2n}$ is dense in
$\mathbb{R}^{2n}$. Given bracket generating property of the pair
(\ref{pervef}) for a generic potential $\Phi$, we conclude according
to Remark~\ref{dense} that this attainable set  coincides with
$\mathbb{R}^{2n}$.

\subsection{Proof of the Proposition~\ref{bezpm}}
\label{propost}

First note that all  points of $\Pi$  are {\it nonwandering} for the
vector field $f$, defined by (\ref{pervef}). Recall that a point $x
\in \mathbb{R}^{2n}$ is nonwandering for $f$ (see \cite[\S6.2]{AM})
if for each neighborhood $U \supset x$ and each $T>0$ there exists
$t>T$ such that $e^{tf}(U)\bigcap U \neq \emptyset$. We will prove
in a moment (Lemma~\ref{nonw}).

Basing on this property we conclude that for each point $x \in \Pi$
and any $t>0$ the points $e^{-tf}(x)$ (contained in $\Pi$)  are
arbitrarily well approximable by points $e^{\tau f}$ with $\tau
>0$.

Acting by a composition of diffeomorphisms $P_N \circ \cdots \circ
P_1 $ belonging to the family (\ref{trifam}) on a point $\tilde{x}
\in \Pi$ we pick the factors   $P_i=e^{-t_if} \ (t>0)$. Each
diffeomorphism  $P_i$ is applied to  a point $y_i=(P_{i-1} \circ
\cdots \circ P_1)(\tilde{x})$ which belongs to $\Pi$.  By
nonwandering property in $\Pi$  we can approximate the action of
$P_i$ by an action of some diffeomorphism $\hat{P}_i=e^{\theta_i f},
\ \theta
>0$.

Thus we proved that  positive orbit of the family (\ref{corfam}) is
dense in $\Pi$ (which is positive orbit of the family (\ref{trifam})
and hence coincides with $\Pi$ given the fact that
$\{\Ad\left(e^{g}\right)f,f\}$, restricted to $\Pi$, form a bracket
generating pair of vector fields on $\Pi$.

\begin{lemma}
\label{nonw} Each point of the hyperplane $\Pi$ is nonwandering for
the vector field $f. \ \square$
\end{lemma}

We will derive this property from Poincare theorem
(\cite[\S3.4]{AM}). Indeed the hyperplane $\Pi$ of zero momentum is
invariant for the Hamiltonian vector field $f$; according to
\cite[\S3.4]{AM} one can introduce a volume form on $\Pi$, which is
preserved by the flow of $f$.

Let us introduce the planes $\Pi_Q=\{\sum_{i=1}^np_i=0, \
\sum_{i=1}^nq_i=Q\}$ and consider the Lebesgue sets $\{H^p \leq c\}$
of the Hamiltonian (\ref{ham-fi}). We will prove in a second
(Lemma~\ref{lecomp}) that intersections of the Lebesgue sets with
each $\Pi_Q$ are compact.

Taking this for granted we see that for each $a,c>0$ the sets
$$\bigcup_{|Q|\leq a}\Pi_Q \bigcap \{H^p \leq c\}$$ are compact and
invariant with respect to the volume-preserving (and
Ha\-mil\-to\-nian-preserving) flow of the vector field $f$. We are
under conditions of Poincare theorem according to which $\forall
Q,c$ points of $\Pi_Q \bigcap \{H^p \leq c\}$ are non-wandering. It
rests to note that each point of $\Pi$ is included in some set
$\Pi_Q \bigcap \{H^p \leq c\}$.

\begin{lemma}
\label{lecomp} Intersections of the Lebesgue sets of the Hamiltonian
$H^p$ with the planes $\Pi_Q=\{\sum_{i=1}^np_i=0, \
\sum_{i=1}^nq_i=Q\}$ are compact.
\end{lemma}

{\it Proof.} Closedness  of the Lebesgue sets $\{H^p \leq c\}$ is
obvious; we prove their boundedness.

Since $\sum_{j=1}^{n-1}\Phi(q_j-q_{j+1})+\Phi(q_n-q_1)$ is bounded
below, say by $-B \leq 0$,
   then the inequality $H^p \leq c$ implies the constraints:
\begin{equation}
  \|p\|^2 \leq c+B , \
\sum_{j=1}^{n-1}\Phi(q_j-q_{j+1})+\Phi(q_n-q_1)
    \leq c.\label{fibo}
\end{equation}
By lower boundedness of the function $\Phi$ and by the growth
conditions  (\ref{grof}) we derive from the second one of the
relations (\ref{fibo})
\begin{equation}
   q_1 -q_2 \leq b
   \bigwedge \cdots q_{n-1}-q_n \leq b \bigwedge q_n-q_1 \leq b,
\label{qb}
\end{equation}
for some constant $b$.

Summing the first $k$ inequalities at the right-hand side of the
implication (\ref{qb}) we conclude
\begin{equation}\label{q1up}
q_1 \leq q_k + (k-1)b, \ k=1, \ldots , n,
\end{equation}
 while summing $n+1-k$
inequalities, starting from the last one, we obtain
\begin{equation}\label{q1lo}
q_k - (n+1-k)b \leq q_1, \ k=1, \ldots , n.
\end{equation}

If we restrict our consideration onto the plane $\Pi_Q$ and sum
separately the inequalities (\ref{q1up}) and (\ref{q1lo}) we get
$$n q_1 \leq Q+b(n-1)n/2, \  nq_1 \geq Q-b(n+1)n/2.$$
Due to invariance with respect to the permutations of particles we
conclude
$$ n^{-1}Q-b(n+1)/2 \leq q_j \leq n^{-1}Q+b(n-1)/2,$$
for each coordinate $q_j$ of a point $(q,p) \in \Pi \bigcap \{H^p
\leq c\}. \ \square$

% ----------------------------------------------------------------


\begin{thebibliography}{KKKK}
\bibitem{AM} Abraham R., Marsden J.E., {\em Foundations of
Mechanics}, 2nd Edition, Perseus Books, 1978.
\bibitem{ASkv} Agrachev A.A., Yu.L.Sachkov, {\em Lectures on
Geometric Control Theory}, Springer-Verlag, 2004.
\bibitem{Bo} Bonnard B. Contr\^olabilit\'e des syst\`emes non lin\'eaires. (French)
C. R. Acad. Sci. Paris, S\'er. I Math.  292(1981), no. 10, 535--537.
\bibitem{Jur}  Jurdjevic V. Geometric Control Theory. Cambridge
University Press, 1997.
\bibitem{Mos} Moser J. Finitely many mass points on the line under
in fluence of an exponential potential -- an integrable system.
Dynamical Systems, Theory and Applications: Battelle Seattle 1974
Rencontres. Editor: J. Moser, Lecture Notes in Physics, vol. 38,
p.467-497.
\bibitem{Pol} Polyanin A.D., Zaitsev V.F. Exact solutions for ordinary differential equations,
2nd edition, Chapman and Hall, 2002.
\bibitem{Sar} Sarychev A.V. Controlling Multiparticle System on the Line. I
(to appear in this Journal).
\end{thebibliography}
\end{document}